\documentclass[12pt,twoside]{article}

\usepackage{a4}
\usepackage{amssymb}
\usepackage{isolatin1}
\usepackage{epsfig}
\usepackage{color}

\newtheorem{thm}{Theorem}
\newtheorem{example}[thm]{Example}

\newenvironment{exam}
        {\begin{example} \quad \rm}
        {\end{example}}

\newenvironment{proof}
        {\pagebreak[2] \vspace{-1pt}{\bf Proof.}  }
        {\hfill $\blacksquare$\vspace{2pt}}

\newenvironment{proofof}[1]
        {\pagebreak[2] \vspace{-1pt}{\bf Proof of #1.}  }
        {\hfill $\blacksquare$ \vspace{2pt}}

\def\nat{{\rm I\! N}}

\def\re{{\mathbb R}}
\def\rat{{\mathbb Q}}

\def\l{\left}
\def\r{\right}
\def\gl{\left\{}
\def\gr{\right\}}
\def\kl{\left(}
\def\kr{\right)}

\def\limk{\lim_{k\to\infty}}
\def\limn{\lim_{n\to\infty}}

\def\In{\subseteq}

\def\mi{\setminus}
\def\abb{\longrightarrow}

\renewcommand{\rho}{\varrho}
\renewcommand{\phi}{\varphi}
\renewcommand{\epsilon}{\varepsilon}

\def\beq{\begin{equation}}
\def\eeq{\end{equation}}
\def\beqar{\begin{eqnarray}}
\def\eeqar{\end{eqnarray}}
\def\beqaro{\begin{eqnarray*}}
\def\eeqaro{\end{eqnarray*}}
\def\bsat{\begin{thm}}
\def\esat{\end{thm}}
\def\blem{\begin{lem}}
\def\elem{\end{lem}}
\def\bkor{\begin{corollary}}
\def\ekor{\end{corollary}}
\def\bprop{\begin{proposition}}
\def\eprop{\end{proposition}}
\def\bdefin{\begin{definition}}
\def\edefin{\end{definition}}
\def\bbew{\begin{proof}}
\def\ebew{\end{proof}}
\def\bbewo{\begin{proofof}}
\def\ebewo{\end{proofof}}
\def\bex{\begin{exam}}
\def\eex{\end{exam}}
\def\bgar{\begin{array}}
\def\ear{\end{array}}

\renewcommand{\rho}{\varrho}
\renewcommand{\phi}{\varphi}

\begin{document}

\title{Some counterexamples on the behaviour of real-valued
  functions and their derivatives}

\author{J\"urgen Grahl and Shahar Nevo}

\date{2016, August 29}

\maketitle

\begin{abstract}
  We discuss some surprising phenomena from basic calculus related to
  oscillating functions and to the theorem on the differentiability of
  inverse functions. Among other things, we see that a continuously
  differentiable function with a strict minimum doesn't have to be
  decreasing to the left nor increasing to the right of the minimum,
  we present a function $f$ whose derivative is discontinuous at one
  point $x_0$ and oscillates only above $f'(x_0)$ (i.e. $f'$ has a
  strict minimum at $x_0$), we compare several definitions of
  inflection point, and we discuss a general version of the theorem on
  the derivative of inverse functions where continuity of the inverse
  function is assumed merely at one point. 
\end{abstract}


\vspace{6pt}

When preparing a lecture on basic calculus you are sometimes reminded
of Hamlet's famous words {\sl ``There are more things in heaven and
  earth, Horatio, than are dreamt of in your philosophy'':} You are
aware of the fact that differentiable real functions are not as tame
as they might appear at first sight and that there are more surprising
phenomena and pitfalls related to them than you would ever have
imagined as a student, and you are faced with the challenge of finding
illustrative examples that give your students a good impression of
those pitfalls.  In this paper we'd like to draw the
  reader's attention to some surprising observations on differentiable
  functions of one variable which have arisen from our own experience
  as lecturers.  Some of them seem to be little-known or even unknown,
  others can already be found in at least one of the many excellent
  calculus textbooks but are included since we consider them
  particularly instructive. We hope that this note is of interest to
  those who want to flavour their lectures with some unexpected or
  even weird examples. In this context we also refer the reader to the
  books of Appell \cite{AppellA}, Gelbaum/Olmsted
  \cite{GelbaumOlmsted} and Rajwade/Bhandari \cite{Rajwade} which
  contain rich collections of counterexamples and pecularities from
  many fields of real analysis. 

\section{Some more things you can do with $x\mapsto x^k\sin\frac{1}{x^m}$}

\label{sec:Oscill}

The functions $x\mapsto x^k\sin\frac{1}{x^m}$ (where $k\ge 0, m\ge 1$)
and related functions are a popular source of enlightening
counterexamples in basic calculus, see, for example,
\cite[Chapter~3]{GelbaumOlmsted}) and for a more detailed discussion
of the properties of those functions \cite[Section 4.5]{AppellA},
\cite{AppellB} and \cite{AppellRoos}.  For instance, the function
$x\mapsto x^2\sin\frac{1}{x^2}$ shows that derivatives are not
necessarily continuous, that they can even be unbounded on compact
intervals, and that differentiable functions need not be rectifiable
on compact intervals, i.e. their graphs can have infinite length.
(From the point of view of complex analysis, the interesting
properties of these functions are of course related to the fact that
$z\mapsto z^k \sin\frac{1}{z^m}$ has an essential singularity at
$z=0$.)

Here we'd like to present some further applications of functions of
this kind. 

\bex  \label{Ex:AccumCritical}
Expressing monotonicity in terms of derivatives is very easy and
elegant: A differentiable function $f:I\abb\re$ on an interval $I$ is
non-decreasing if and only if $f'(x)\ge 0$ for all $x\in I$.
Things are a bit less elegant for {\it strict} monotonicity: The
condition $f'(x)>0$ for all $x\in I$ is sufficient but not necessary
for $f$ being strictly increasing. In fact, as known from elementary
calculus, $f$ is strictly increasing if and only if $f'(x)\ge 0$ for
all $x\in I$ and if there doesn't exist a proper interval $J\In I$
such that $f'(x)=0$ for all $x\in J$. 

\begin{figure}[htb]
  \begin{center}
\includegraphics[width=0.475\textwidth]{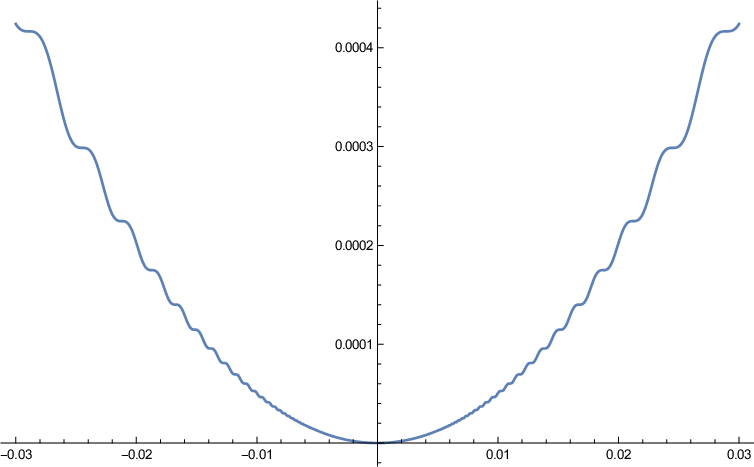}\hfill
\includegraphics[width=0.475\textwidth]{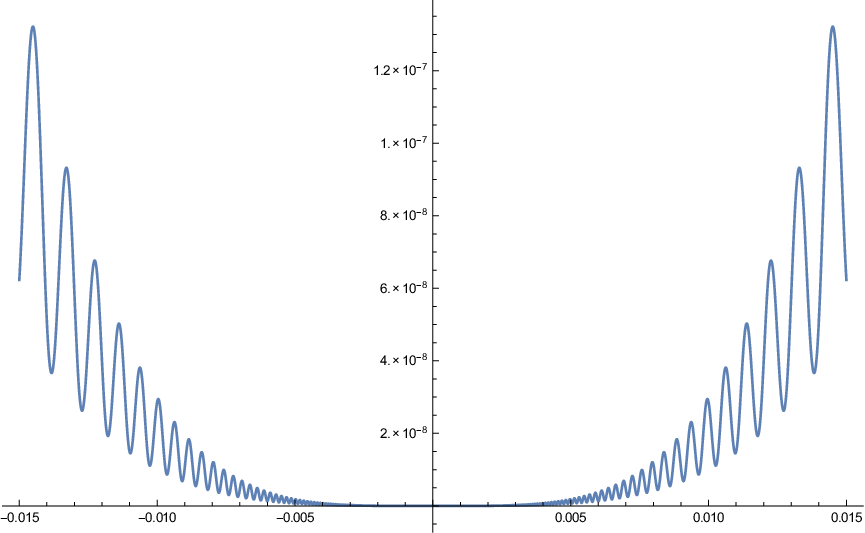}
    \caption{{\sf (a) $x\mapsto\int_{0}^{x} t\cdot \kl
        1+\cos\frac{1}{t}\kr\,dt$ \qquad (b) $x\mapsto x^4\kl2+\sin\frac{1}{x}\kr$}}  
    \label{fig:AccumCritical}    \label{fig:MinNichtMonoton}
  \end{center}
\end{figure}

In particular, there must be strictly increasing functions which have
infinitely many critical points within a compact interval. An example
is provided by the function $f:\re\abb\re$ defined by
$$f(x):=\int_{0}^{x} t\cdot \kl 1+\cos\frac{1}{t}\kr\,dt \qquad\mbox{
  for all } x\in\re $$
(Figure \ref{fig:AccumCritical} (a)) which is
obviously differentiable and whose critical points are 0 and the
points $\frac{1}{(2n+1)\pi}$ (where $n\in\nat$) which accumulate at
the origin. 

In \cite[p.  186]{koehler} a similar example is
  given: Define $f : [0, 1] \rightarrow \mathbb{R}$ by $f(0) = 0$ and
$$f(x) := x \cdot \kl 2 - \sin(\log(x)) - \cos(\log(x)) \kr \quad
\mbox{ for } 0 < x \leq 1 .$$
An easy calculation shows that this function is strictly increasing
with critical points at $e^{-2\pi n}$ ($n\in\nat$). However, at $x=0$
it is only continuous, but not differentiable.

It's worth noting that the (almost trivial) monotonicity criterion
mentioned above can be extended as follows: {\sl Let $f:[a,b]\abb\re$
  be a continuous function. Assume that there is a countable set $S$
  such that $f$ is differentiable on $[a,b]\mi S$ and $f'(x)>0$ for
  all $x\in [a,b]\mi S$. Then $f$ is strictly monotonically increasing
  on $[a,b]$.} An amazingly short proof for this classical result was
given by L.~Zalcman \cite{Zalcman}.
\eex

\bex \label{Ex:MinNonMon}
If $f:\re\abb\re$ is a continuously differentiable function which has
a strict local minimum at some point $x_0$, then most students tend to
expect that there is a small neighbourhood $(x_0-\delta,x_0+\delta)$ such that $f$ is
decreasing on $(x_0-\delta,x_0]$ and increasing on
$[x_0,x_0+\delta)$. 

The function \cite[p.~36]{GelbaumOlmsted} 
$$f(x):=\left\{\begin{array}{ll} x^4\kl2+\sin\frac{1}{x}\kr & \mbox{ for }
  x\in\re\mi\gl 0\gr, \\[9pt] 0 & \mbox{ for }
  x=0\end{array}\right. $$
shows that this is wrong (Figure \ref{fig:MinNichtMonoton}
  (b)). Obviously, $f$ has a strict local minimum at 0, and $f$ is
  even continuously differentiable on $\re$ with $f'(0)=0$, but $f'$
  assumes both positive and negative values in any neighbourhood of 0
  as we can easily see from 
$$f'(x)=x^2\cdot\kl 8x+4x\sin\frac{1}{x}-\cos\frac{1}{x}\kr\qquad \mbox{ for }
x\ne0.$$
\eex

\begin{figure}[htb]
  \begin{center}
\includegraphics[width=0.475\textwidth]{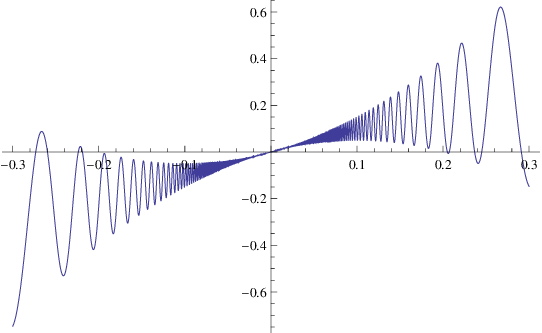}\hfill
\includegraphics[width=0.475\textwidth]{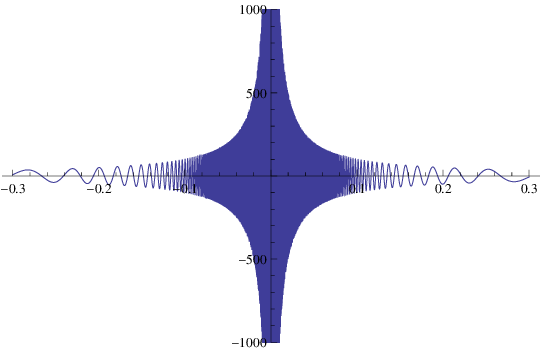}
    \caption{{\sf $x\mapsto x+5x^2\sin\frac{1}{x^2}$ and its derivative}}  
    \label{fig:NotInjective}
  \end{center}
\end{figure}

\bex\label{Ex:NotInjective}
If $f:\re\abb\re$ is differentiable and $f'(x_0)\ne 0$, this does not
imply that $f$ is monotonic in some neighbourhood of $x_0$. This is
illustrated by the function \cite[p.~37]{GelbaumOlmsted} 
$$f:\re\abb\re, \quad f(x):=\gl\bgar{ll} x+\alpha x^2\sin\frac{1}{x^2}
  & \mbox{ for } x\ne 0\\[2pt] 0 & \mbox{ for } x=0,\ear\r. \qquad\mbox{
  where } \alpha\ne 0$$
which is differentiable everywhere with $f'(0)=1\ne 0$, but is not
monotonic (hence not one-to-one) on any neighbourhood of 0 since $f'$ assumes both positive
  and negative values there. In fact, $f'$ is unbounded in both
  directions on each such neighbourhood (see Figure \ref{fig:NotInjective}). 
\eex

\bex\label{Ex:OscillAbove}
By Darboux's intermediate value theorem for derivatives, if a
derivative $f'$ is not continuous at some point $x_0$, then $x_0$ has
to be an oscillation point of $f'$ (i.e. there is a certain interval
$J$ of positive length such that $f'$ assumes all values in $J$ in
arbitrarily small neighbourhoods of $x_0$). In the usual examples for
this situation (like $x\mapsto x^2\sin\frac{1}{x}$), $f'$ oscillates
in both directions, above and below $f'(x_0)$, i.e. for arbitrary
small neighbourhoods $U$ of $x_0$, $f'(x_0)$ is an inner point of
$f'(U)$.

\begin{figure}[htb]
  \begin{center}
\includegraphics[width=0.475\textwidth]{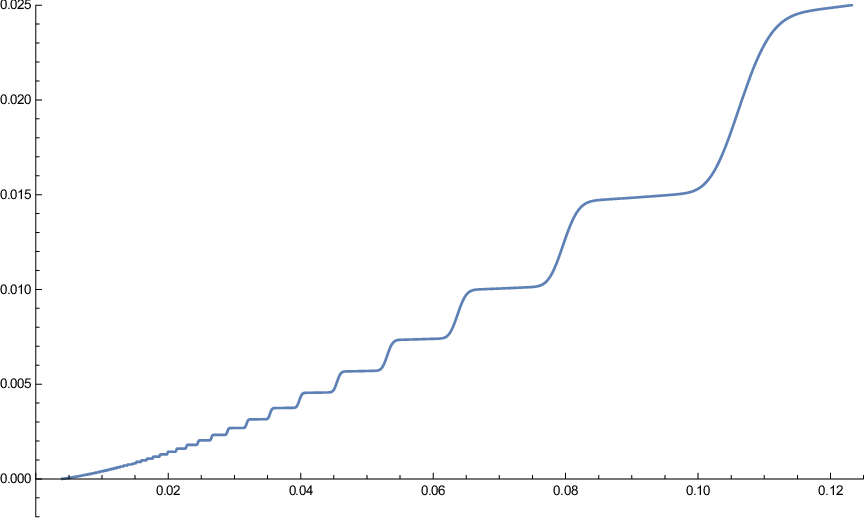}\hfill
\includegraphics[width=0.475\textwidth]{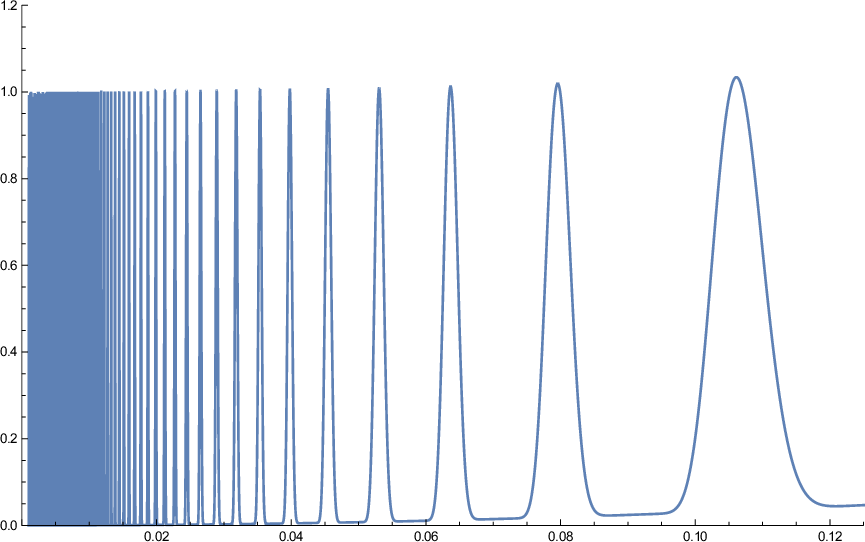}

\caption{{\sf $x \mapsto
    x^3+\int_{0}^{x}\l|\cos\frac{1}{t}\r|^{1/|t|}\,dt$ and its
    derivative (Note the different scales in both figures.)}}  
\label{fig:IntegralCos}
\end{center}
\end{figure}

The function
$$f(x):=x^3+\int_{0}^{x} \l|\cos\frac{1}{t}\r|^{1/|t|}\,dt$$
(a slight variation of the function from Example \ref{Ex:AccumCritical}) 
is an example of a differentiable function $f:\re\abb\re$ which
satisfies $f'(0)=0< f'(x)$ for all $x\in\re\mi\gl0\gr$ and whose derivative
$f'$ is not continuous at 0, i.e.~$f'$ oscillates only above $f'(0)$
(see Figure \ref{fig:IntegralCos}). In other words, $f'$ has even a
strict minimum at the origin. Here, $f'(0)=0$ is not obvious, so we provide the details.  

\bbew{} 
The function $t\mapsto \l|\cos\frac{1}{t}\r|^{1/|t|}$ is bounded and
continuous almost everywhere, so by Lebesgue's criterion it is integrable (even in the
sense of Riemann) on compact intervals. Therefore, $f$ is
well-defined. For $x\ne 0$ it is clear from the fundamental theorem of
calculus that $f'(x)=3x^2+\l|\cos\frac{1}{x}\r|^{1/|x|}> 0$. 

So it remains to show that $f'(0)=0$. For this purpose we substitute
$u:=\frac{1}{t}$ to write $f(x)$ in the form 
$$f(x)=x^3+\int_{1/x}^{\infty} \frac{|\cos u|^u}{u^2}\,du \qquad
\mbox{ for } x>0$$
and consider the points  
$$a_k:= k\pi, \qquad b_k:=\kl k+\frac{1}{k^{1/4}}\kr\cdot \pi, \qquad
c_k:=\kl k+1-\frac{1}{k^{1/4}}\kr\cdot \pi \qquad (k\in\nat).$$
For $k\ge 16$ we have $a_k\le b_k\le c_k\le a_{k+1}$, 
$$\int_{a_k}^{b_k} \frac{|\cos u|^u}{u^2}\,du
\le \frac{b_k-a_k}{a_k^2}=\frac{1}{\pi \cdot k^{9/4}}$$
and
$$\int_{c_k}^{a_{k+1}} \frac{|\cos u|^u}{u^2}\,du
\le\frac{a_{k+1}-c_k}{c_k^2}\le\frac{1}{\pi \cdot k^{9/4}}.$$
Furthermore, using the estimate $\cos y\le 1-\frac{y^2}{4}$ for $0\le
y\le \frac{\pi}{2}$, for $k\ge 16$ we obtain
$$\int_{b_k}^{c_k} \frac{|\cos u|^u}{u^2}\,du
\le \frac{c_k-b_k}{b_k^2}\cdot\kl\cos \frac{\pi}{k^{1/4}}\kr^{b_k}
\le \frac{1}{\pi k^2}\cdot\kl 1-\frac{\pi^2}{4\sqrt{k}}\kr^{k\pi}.$$ 
Here, $\limk \kl 1-\frac{\pi^2}{4\sqrt{k}}\kr^{\sqrt{k}\pi}=e^{-\pi^3/4}<1$,
so if we set $\alpha:=\frac{1}{2}\cdot \kl1+e^{-\pi^3/4}\kr$, then
there exists some $N_0\ge 16$ such that
$$ \kl 1-\frac{\pi^2}{4\sqrt{k}}\kr^{\sqrt{k}\pi}\le \alpha \qquad \mbox{ for all }
k\ge N_0,$$ 
and we get 
$$\int_{b_k}^{c_k} \frac{|\cos u|^u}{u^2}\,du 
\le \frac{1}{\pi k^2}\cdot \alpha^{\sqrt{k}} 
\qquad \mbox{ for all } k\ge N_0.$$ 
Combining these estimates, for all $N\ge N_0-1$ we deduce
\beqaro
0&\le& \int_{(N+1)\pi}^{\infty} \frac{|\cos u|^u}{u^2}\,du \\
&=&\sum_{k=N+1}^{\infty} \kl\int_{a_k}^{b_k} \frac{|\cos u|^u}{u^2}\,du 
+\int_{b_k}^{c_k} \frac{|\cos u|^u}{u^2}\,du 
+\int_{c_k}^{a_{k+1}} \frac{|\cos u|^u}{u^2}\,du \kr\\
&\le & \sum_{k=N+1}^{\infty}\kl \frac{2}{\pi \cdot k^{9/4}}
+\frac{1}{\pi k^2}\cdot \alpha^{\sqrt{k}} \kr
\le  \frac{2}{\pi} \int_{N}^{\infty} \frac{dx}{x^{9/4}}
+ \frac{1}{\pi}  \int_{N}^{\infty} \frac{\alpha^{\sqrt{x}}}{\sqrt{x}}\; dx\\
&=& \frac{8}{5\pi}\cdot\frac{1}{N^{5/4}}+\frac{2}{\pi}
\int_{\sqrt{N}}^{\infty} \alpha{^y}\,dy
=\frac{8}{5\pi}\cdot\frac{1}{N^{5/4}}+\frac{2\alpha^{\sqrt{N}}}{\pi\log\frac{1}{\alpha}}.
\eeqaro
Now let some $x>0$ be given. W.l.o.g. we may assume $x<\frac{1}{N_0
  \pi}$. Then there exists some $N_x\in\nat$, $N_x\ge N_0$ such that
$\frac{1}{(N_x+1)\pi} <x \le \frac{1}{N_x\pi}$. Here for $x\to 0+$ we
have $N_x\abb \infty$, so we obtain
\beqaro
0&\le& \frac{f(x)}{x}
\le x^2+(N_x+1)\pi\cdot\int_{N_x\pi}^\infty
\frac{\l|\cos u\r|^{u}}{u^2}\,du\\
&\le&x^2+\frac{8}{5}\cdot\frac{N_x+1}{(N_x-1)^{5/4}}+\frac{2}{\log\frac{1}{\alpha}}\cdot
(N_x+1)\cdot\alpha^{\sqrt{N_x-1}} \quad\abb 0 \mbox{ for } x\to0+.
\eeqaro
This shows that $\frac{f(x)}{x}\abb 0$ for $x\to0+$. In view of
$f(-x)=-f(x)$ the same holds for $x\to0-$. So $f$ is differentiable at
0 with $f'(0)=0$.  

It is clear that $f'$ is not continuous at 0. 
\ebew

This example can be easily modified such that $f'$ is even unbounded
near the origin (and still oscillates only in the positive direction, 
i.e.~between 0 and $\infty$), by setting 
$$f(x):=\int_{0}^{x}
\frac{1}{\sqrt{|t|}}\cdot\l|\cos\frac{1}{t}\r|^{1/|t|^2}\,dt.$$
Furthermore, functions with the properties discussed in this item can
also be constructed as integral functions of piecewise linear
functions, having increasingly thinner and increasingly higher peaks
which accumulate at the origin but this might be considered less
elegant than the examples discussed above.
\eex

\bex\label{Ex:Inflection}
In calculus textbooks, there are various definitions of an {\it
  inflection point}\footnote{For a more detailed discussion of five
  different definitions of inflection point we refer to \cite[Section
  5.5]{Rajwade}. Here we just want to point out some noteworthy
  aspects.}. Maybe the most common one (see, for example \cite[p.
195]{koehler}, \cite[p. 147]{koenigsberger} and \cite[p.
186]{SalasHille}) is the following one.

\begin{itemize}
\item[(D1)]  
{\sl A continuous function $f:[a,b]\abb\re$ has an inflection point at
$x_0\in(a,b)$ if $f$ is strictly convex on one side of $x_0$ (more
precisely, in a certain interval $(x_0-\delta,x_0)$ resp. $(x_0,
x_0+\delta)$) and strictly concave on the other side.} 
\end{itemize}
According to this definition also a function like
\beq\label{inflection}
f(x):=\gl\bgar{ll} x^2 & \mbox{ for } x<0 \\ \sqrt{x} & \mbox{ for }
x\ge 0\ear\r.
\eeq
has an inflection point at $x=0$ though it is not differentiable
there. 

Sometimes (for example in \cite[p. 388]{Apostol})
also the following definition can be found:  
\begin{itemize}
\item[(D2)]  
{\sl A function $f:[a,b]\abb\re$ has an inflection  point at
  $x_0\in(a,b)$ if $f$ is differentiable at $x_0$ and if the 
  graph of $f$ is strictly above the tangent line $x\mapsto
  f(x_0)+f'(x_0)\cdot(x-x_0)$ on one side of $x_0$ and strictly below this
tangent line on the other side}. 
\end{itemize}

\begin{figure}[htb]
  \begin{center}
\includegraphics[width=0.7\textwidth]{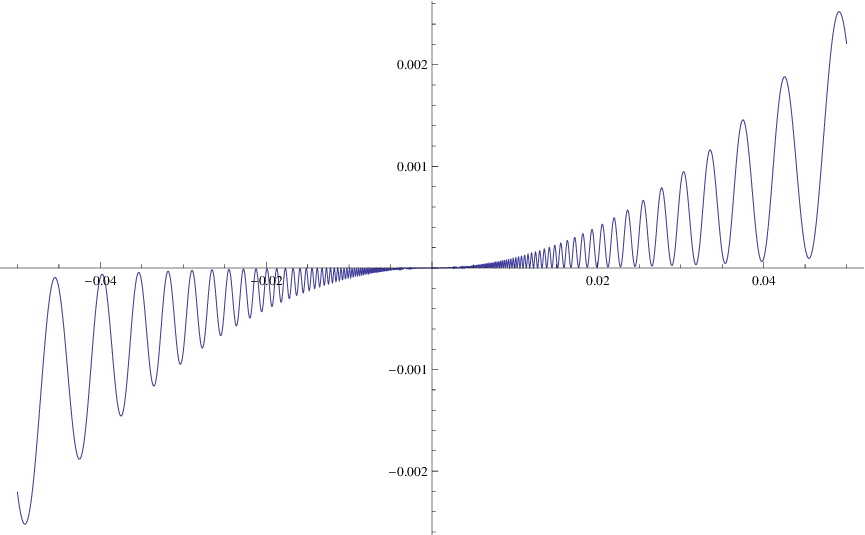}
    \caption{{\sf $x\mapsto x^3+{\rm sgn}(x)\cdot x^2\sin^2\frac{1}{x}$}}  
    \label{fig:InflectionPoints}
  \end{center}
\end{figure}

This definition does not apply to the function in (\ref{inflection}),
but in the case of functions which are differentiable at $x_0$, it is
more general than the first definition. For example, if we set
$$f(x):=
x^3+{\rm sgn}(x)\cdot x^2\sin^2\frac{1}{x}
=\gl\bgar{ll} x^3+x^2\sin^2\frac{1}{x} & \mbox{ for } x>0, \\ 
0 & \mbox{ for } x=0, \\
x^3-x^2\sin^2\frac{1}{x} & \mbox{ for } x<0, \ear\r.$$
then $f$ is differentiable at $x=0$ with $f'(0)=0$, $f(x)>0$ for all 
$x>0$ and $f(x)<0$ for all $x<0$, so $f$ has an inflection point at $x=0$
in the sense of (D2). However, $f$ is neither convex nor concave in
$(0,\delta)$ for any $\delta>0$, and the same holds in the intervals
$(-\delta,0)$ (see Figure \ref{fig:InflectionPoints}). 

Definition (D2) even applies to functions which are not continuous in
  any neighbourhood of the inflection 
  point. For example, with this definition the function 
$$f(x):=\gl\bgar{ll} x^3 & \mbox{ for } x\in\rat, \\ 
x^5 & \mbox{ for } x\in\re\mi\rat.\ear\r.$$
(whose graph looks like the union of the graphs of $x\mapsto x^3$ and
of $x\mapsto x^5$, though this is of course a visual illusion) has
an inflection point at $x=0$.

Finally, some authors (for example \cite[p. 150]{Erwe}) prefer the
following definition which is, however, the least general one since it
is restricted to functions that are differentiable everywhere (not
just at the inflection point):  
\begin{itemize}
\item[(D3)] 
{\sl A differentiable function $f:[a,b]\abb\re$ has an inflection point at
  $x_0\in\;(a,b)$ if $f'$ has a strict local extremum at $x_0$.} 
\end{itemize}

For functions that are differentiable everywhere, definition (D3) is
more general than (D1) but less general than (D2). More precisely, if
$f:[a,b]\abb\re$ is differentiable, then 
\begin{itemize}
\item 
if $x_0\in\;(a,b)$ is an inflection point in the sense of (D1), it is
also an inflection point in the sense of (D3) (since strict convexity
resp. concavity can be described in terms of the first derivative
being strictly increasing resp. decreasing), 
\item 
and if it is an inflection point in the sense of (D3), it is also an
inflection point in the sense of (D2) as easily seen by considering
the function 
$$d(x):=f(x)-f(x_0)-f'(x_0)(x-x_0)$$
which under the conditions of (D2) is strictly monotonic with a zero at
$x_0$, hence changes its sign at $x=x_0$. 
\end{itemize}

\begin{figure}[htb]
  \begin{center}
\includegraphics[width=0.7\textwidth]{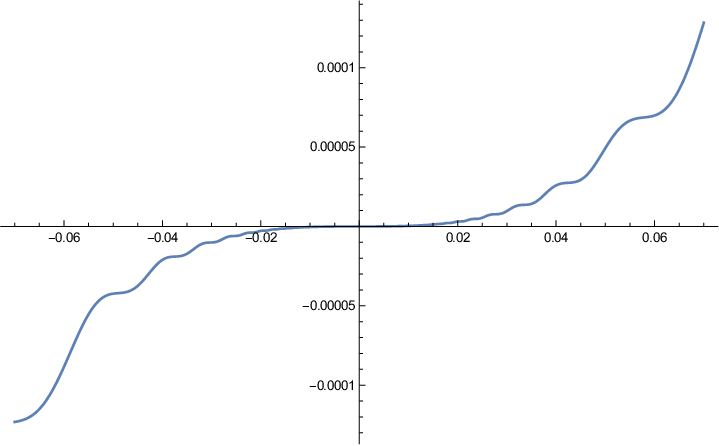}
    \caption{{\sf $x\mapsto\int_{0}^{x}  t^2\kl1.1+\sin\frac{1}{t}\kr \,dt$}}  
    \label{fig:InflectionP-NonMonoton}
  \end{center}
\end{figure}

That an inflection point in the sense of (D2) need not be an
inflection point in the sense of (D3) is again illustrated by the
function from Figure \ref{fig:InflectionPoints}. To find an example
for an inflection point in the sense of (D3), but not in the sense of
(D1), we recall that a strict local minimum doesn't necessarily mean a
change from monotonically decreasing to monotonically increasing as we
know from Example \ref{Ex:MinNonMon}. In fact, the desired example is
provided by an antiderivative of a function as in Example
\ref{Ex:MinNonMon}, say the function $f:\re\abb\re$ defined by 
$$f(x):=\int_{0}^{x} t^2\kl1.1+\sin\frac{1}{t}\kr \,dt \qquad\mbox{
  for } x\in\re.$$
(Here we have slightly modified the function from
Example \ref{Ex:MinNonMon} in order to make the desired phenomenon
more apparent in the graph. But it still has the very same properties:
It has a strict extremum at the origin, but isn't decreasing on one
side nor decreasing on the other side of the extremum.)

Then $f$ is differentiable with $f'(x)>f'(0)=0$ for all $x\ne0$, so
(D3) is satisfied for $x_0=0$. However, (D1) is violated since $f''$
changes its sign in arbitrary small neighborhoods of $x_0$ (see Figure
\ref{fig:InflectionP-NonMonoton}).  
\eex

\bex\label{Ex:IntegrandUnbounded}
There are certain analogies between infinite series and improper
Riemann integrals of the form $\int_{0}^{\infty} f(t)\,dt$. For
example, if $f:[0,\infty[\abb[0,\infty[$ is monotonically decreasing,
then the convergence of the improper integral $\int_{0}^{\infty}
f(t)\,dt$ is equivalent to the convergence of the series
$\sum_{n=0}^{\infty} f(n)$. 

\begin{figure}[htb]
\begin{center}
\includegraphics[width=0.475\textwidth]{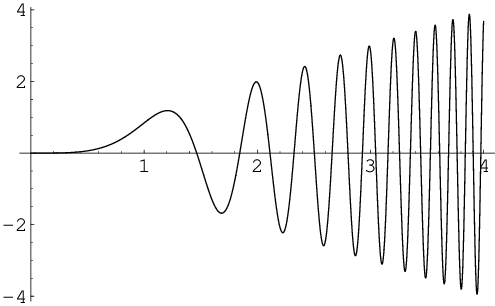}
\hfill
\includegraphics[width=0.475\textwidth]{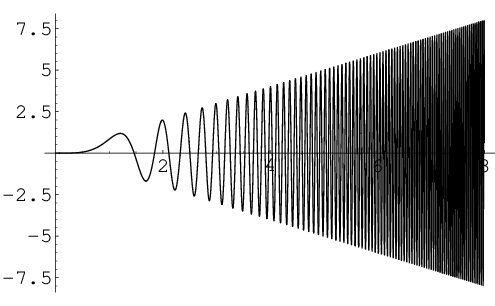} 
\caption{\sf $x\mapsto x\sin (x^3)$ in the intervals $[0;4]$ and $[0;8]$} 
\label{uneigUnbeschr}
\end{center}
\end{figure}

Now if $\sum_{n=0}^{\infty} a_n$ is a convergent series, $(a_n)_n$
necessarily converges to 0. This might (mis)lead to the conjecture
that if the improper Riemann integral $\int_{0}^{\infty} f(t)\,dt$
converges (where $f$ is a continuous function) then $\lim_{t\to\infty}
f(t)=0$. But this is wrong. On the contrary, the function
$$f(x):=x\sin (x^3)$$
(Figure \ref{uneigUnbeschr}) shows that $f(x)$
can be even unbounded as $x\to\infty$ (cf. also \cite[p.
46]{GelbaumOlmsted}). Here, to show the convergence of
$\int_{0}^{\infty} f(t)\,dt$ we substitute $y:=t^3$ which gives
$$\int_{1}^{R} f(t)\,dt
=\int_{1}^{R^3} \frac{\sin y}{3y^{1/3}}\,dy \qquad \mbox{ for all }
 R>0$$
and observe that the improper integral 
$\int_{1}^{\infty} \frac{\sin y}{3y^{1/3}}\,dy $ converges.
\eex

\bex\label{Ex:NoLimitDeriv}
When talking about differentiable functions $f:[a;\infty)\abb\re$ with
a certain asymptotic behaviour near $\infty$, say $\lim_{x\to\infty}
f(x)=0$, the question arises how the derivative $f'$ behaves near
$\infty$. Probably some students will think that $f'(x)$ is bound to
tend to 0 as $x\to\infty$.  

However, there's no good reason to think so, as you can illustrate by
the following analogy: If a small particle is trapped in a small room
and we make this room smaller and smaller, this particle might still
be able to move arbitrarily fast within this room (as long as it is
able to accelerate and decelerate within very short times).

\begin{figure}[htb]
  \begin{center}
\includegraphics[width=0.7\textwidth]{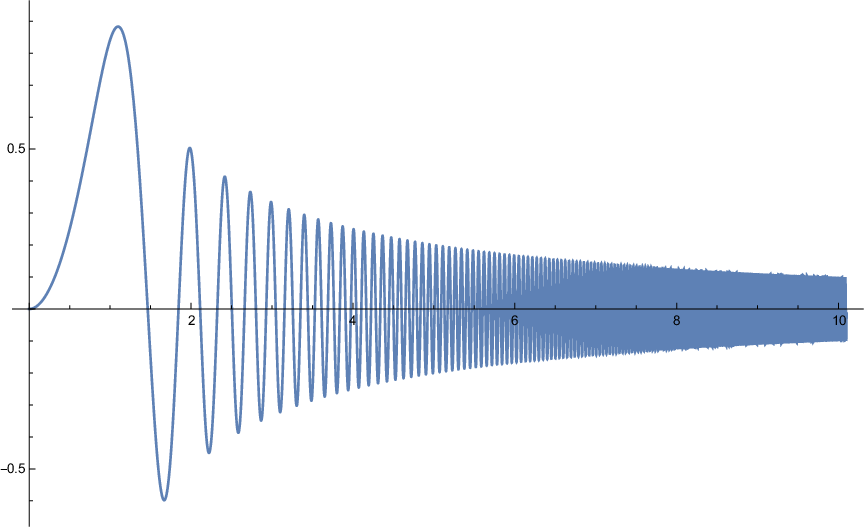}
    \caption{{\sf $x\mapsto \frac{\sin(x^3)}{x}$}}  
    \label{fig:NoLimitDeriv}    
  \end{center}
\end{figure}

In fact, a simple counterexample is provided by the function 
$$f(x):=\frac{\sin(x^3)}{x}$$
which satisfies $\lim_{x\to\infty} f(x)=0$. Here $\lim_{x\to\infty} f'(x)$
does not exist, and $f'$ is even unbounded near $\infty$ (Figure
\ref{fig:NoLimitDeriv}). Actually, $f'$ behaves very much like the
function in Example \ref{Ex:IntegrandUnbounded}.

It is more difficult to give a counterexample where $f$ is even
monotonic, i.e. $f'$ doesn't change signs. Here a near relative of the
function from Example \ref{Ex:MinNonMon} proves to be helpful: Choose
an arbitrary $\varepsilon>0$ and define
$$f(x):=\int_{0}^{x} |\cos t|^{t^{2+\varepsilon}}\,dt \qquad\mbox{ for } x> 0.$$
Then $f$ is differentiable on $(0;\infty)$ with non-negative
derivative $f'(x)=|\cos x|^{x^{2+\varepsilon}}$ (whose zeros are isolated), i.e. $f$
is strictly increasing, and $f'(x)$ does not tend to 0 for
$x\to\infty$. Furthermore, a reasoning similar to Example
\ref{Ex:OscillAbove} shows that the improper integral
$\int_{0}^{\infty} |\cos t|^{t^{2+\varepsilon}}\,dt=\lim_{x\to\infty} f(x)$ exists. 

\begin{figure}[htb]
  \begin{center}
\includegraphics[width=0.475\textwidth]{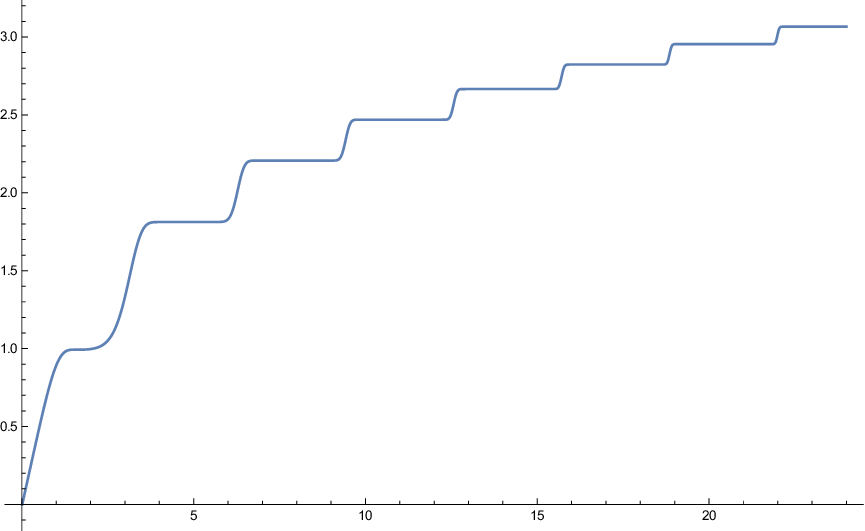}\hfill
\includegraphics[width=0.475\textwidth]{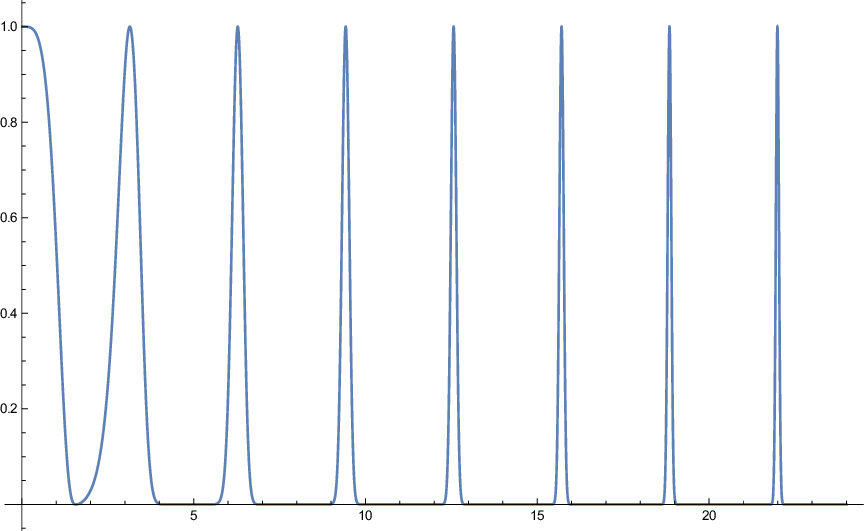}
    \caption{{\sf $x\mapsto \int_{0}^{x} |\cos t|^{t^{2.01}}\,dt$ and
        its derivative}}  
    \label{fig:NoLimitDerivMon}    
  \end{center}
\end{figure}

\bbew
We consider the points    
$$a_k:= k\pi, \qquad b_k:=\kl k+\frac{1}{k^{1+\varepsilon/4}}\kr\cdot \pi, \qquad
c_k:=\kl k+1-\frac{1}{k^{1+\varepsilon/4}}\kr\cdot \pi $$
for $k\in\nat$. 
Then $a_k\le b_k\le c_k\le a_{k+1}$ for $k$ large enough, and we can
estimate 
$$\int_{a_k}^{b_k} |\cos t|^{t^{2+\varepsilon}} \,dt
\le b_k-a_k =\frac{\pi}{k^{1+\varepsilon/4}},\qquad
\int_{c_k}^{a_{k+1}} |\cos t|^{t^{2+\varepsilon}} \,dt
\le \frac{\pi}{k^{1+\varepsilon/4}},$$
$$\int_{b_k}^{c_k} |\cos t|^{t^{2+\varepsilon}} \,dt
\le (c_k-b_k)\cdot\kl\cos \frac{\pi}{k^{1+\varepsilon/4}}\kr^{b_k^{2+\varepsilon}}
\le \pi \cdot\kl 1-\frac{\pi^2}{4k^{2+\varepsilon/2}}\kr^{(k\pi)^{2+\varepsilon}}.$$ 
From $\limk \kl1-\frac{\pi^2}{4k^{2+\varepsilon/2}}\kr^{k^{2+\varepsilon/2}}=e^{-\pi^2/4}<1$ 
we see that there is a constant $\alpha<1$ such that
$$\int_{b_k}^{c_k} |\cos t|^{t^{2+\varepsilon}} \,dt
\le \pi\cdot\alpha^{k^{\varepsilon/2}}$$
for $k$ sufficiently large. Here the series
$\sum_{k=1}^{\infty}\frac{\pi}{k^{1+\varepsilon/4}}$ is obviously
convergent, and the convergence of
$\sum_{k=1}^{\infty}\alpha^{k^{\varepsilon/2}}$ follows from the
convergence of the improper integral $\int_1^\infty
\alpha^{x^{\varepsilon/2}}\,dx=\frac{2}{\varepsilon}\int_1^\infty
\alpha^y y^{2/\varepsilon-1}\,dy$. 
Summarizing, we conclude that $\int_{0}^{\infty} |\cos t|^{t^{2+\varepsilon}}\,dt$
exists. 
\ebew

The graph of $f$ and $f'$ (for $\varepsilon=\frac{1}{100}$) is shown in Figure
\ref{fig:NoLimitDerivMon}.   
\eex

\section{On the differentiability of inverse functions}

The theorem on the inverse function of a differentiable function can
be stated in the following general form.

\bsat{}\label{Thm:DerivInverse}
Let $I\In\re$ be some open interval, $x_0\in I$ and $f: I \abb \re$ be a
function such that 
\begin{itemize}
\item[(i)] $f$ is differentiable at $x_0$ with $f'(x_0)\ne 0$ and
\item[(ii)] $f$ is one-to-one and the inverse $f^{-1}:f(I)\abb I$ is
  continuous at $y_0:=f(x_0)$.
\end{itemize}
Then $f^{-1}$ is differentiable at $y_0$, and its derivative is
$$(f^{-1})'(y_0) = \frac{1}{f'(f^{-1}(y_0))} = \frac{1}{f'(x_0)} .$$
\esat

Of course, here the general definition of derivative is used which
doesn't assume that $f^{-1}:f(I)\abb I$ is defined in a neighborhood
of $y_0$; it just requires that $y_0$ is an accumulation point of
$f(I)$ (which is inevitable since otherwise it wouldn't make sense to
talk about differentiability of $f^{-1}$ at $y_0$).

It seems that this general version of Theorem \ref{Thm:DerivInverse}
is little known; in fact, the only references we could find are in
Hebrew \cite[p. 163-164]{Meizler} and in German \cite[p. 292]{Deiser}.

For the convenience of the reader we recall the short proof of Theorem
\ref{Thm:DerivInverse}.

\bbew{}
First we have to make sure that $y_0=f(x_0)$ is an accumulation point of 
$f(I)$. Indeed, there is a sequence $\gl x_n\gr_n$ in $I\mi\gl x_0\gr$ that tends to $x_0$, and by the
continuity of $f$ at $x_0$ we obtain $f(x_n)\abb f(x_0)$
($n\to\infty$). Since $f$ is one-to-one, none of the values $f(x_n)$
is equal to $f(x_0)$, so $\gl f(x_n)\gr_n$ is a sequence in
$f(I)\mi\gl y_0\gr$ that tends to $y_0$, hence confirming
$y_0$ as an accumulation point of $f(I)$. 

Now we consider an arbitrary sequence $\gl y_n\gr_n$ in $f(I)\mi\gl
y_0\gr$ that tends to $y_0$. Then by the continuity of $f^{-1}$ at $y_0$
we have $x_n:=f^{-1}(y_n)\abb f^{-1}(y_0)=x_0$ for $n\to\infty$. From
this and the existence of $f'(x_0)\ne0$ we obtain
$$\frac{f^{-1}(y_n)-f^{-1}(y_0)}{y_n-y_0}
=\frac{x_n-x_0}{f(x_n)-f(x_0)}\abb \frac{1}{f'(x_0)} \qquad
(n\to\infty).$$ 
Since this holds for any sequence $\gl y_n\gr_n$ in $f(I)\mi\gl
y_0\gr$ with  $\limn y_n=y_0$, we deduce the existence of the
derivative $(f^{-1})'(y_0)=\frac{1}{f'(x_0)}.$
\ebew

Most textbooks on calculus (see, for example
\cite[Theorem~29.9]{Ross}) prefer to give a less general version of
this theorem where condition (ii) is replaced by the following
condition.

\begin{itemize}
\item[(ii)'] $f$ is one-to one and continuous on $I$ (i.e. it is
  strictly monotonic). 
\end{itemize}
Of course, (ii)' implies (ii) since the inverse function of a
continuous strictly monotonic function is continuous.  

The condition that $f$ is one-to-one in (ii) resp.~(ii)' cannot be
skipped since $f'(x_0)\ne 0$ does not imply that $f$ is
monotonic in some neighbourhood of $x_0$ (even if $f$ is
differentiable on the whole of $I$) as we have seen in Example
\ref{Ex:MinNonMon}. 

However, if in Theorem \ref{Thm:DerivInverse} one assumes that $f$ is
differentiable on $I$ with $f'(x)\ne 0$ for all $x\in I$ (not just for
$x=x_0$), then $f$ is one-to-one by Darboux's intermediate value
theorem for derivatives.   

On the other hand, it does not suffice just to assume (in (ii) resp.
(ii)') that $f$ is one-to-one (as it is done, for example, in
\cite{Ball} and \cite{Moise}); Theorem \ref{Thm:DerivInverse} might
fail if $f^{-1}$ is not continuous at $y_0$. This is illustrated by
the following example. 

\bex{} \label{CounterexInverse}
We set $I:=[-2,2]$, $B_n:=\l[\frac{1}{n+1},\frac{2n+1}{2n(n+1)}\kr$
and 
$$A:=[-2,2]\mi\kl \gl0\gr\cup\bigcup_{n=1}^{\infty} B_n
\cup\bigcup_{n=1}^{\infty} (-B_n)\kr.$$
(Observe that $\frac{2n+1}{2n(n+1)}$ is just the middle of the
interval $\l[\frac{1}{n+1},\frac{1}{n}\r]$.) 

\begin{figure}[htb]
  \begin{center}
\includegraphics[width=0.7\textwidth]{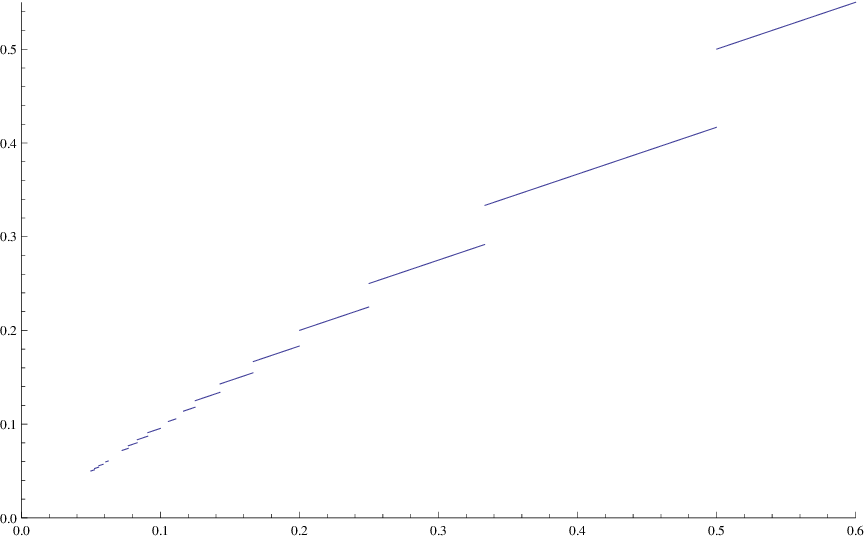}
    \caption{{\sf The graph of the function in Example
        \ref{CounterexInverse} in the interval $[0.05,0.6]$}}  
    \label{fig:CounterexInverse}
  \end{center}
\end{figure}

Then $A$ is an uncountable subset of $\re$, so by the
Cantor-Bernstein-Schroeder theorem \cite{Hinkis} there exists a one-to-one map $T$
of $[-2,-1]\;\cup\;[1,2]$ onto $A$. We define $f:[-2,2]\abb\re$ by
$$f(x):=\gl\bgar{ll}0 & \mbox{ for } x=0, \\[6pt] 
\frac{1}{n+1}+\frac{1}{2}\kl x-\frac{1}{n+1}\kr & \mbox{ for }
\frac{1}{n+1}\le x<\frac{1}{n},\\[6pt]
-\frac{1}{n+1}+\frac{1}{2}\kl x+\frac{1}{n+1}\kr & \mbox{ for }
-\frac{1}{n}< x\le-\frac{1}{n+1},\\[6pt]
T(x)&\mbox{ for } x\in [-2,-1]\;\cup\;[1,2].\ear\r.$$
The graph of $f$ in the interval $[0.05,0.6]$ is sketched in Figure
\ref{fig:CounterexInverse}. Obviously, $f(-x)=-f(x)$ for
$-1<x<1$. 

Since $f$ maps the intervals $\l[\frac{1}{n+1},\frac{1}{n}\kr$
linearly onto $B_n$, it is clear that $f$ maps $I$ onto itself in a
one-to-one fashion. For any given
$x\in\l[\frac{1}{n+1},\frac{1}{n}\kr$  we can write 
$x=\frac{1}{n+1}+\tau$ with $0\le \tau<\frac{1}{n(n+1)}$, so we obtain
$$\frac{f(x)-f(0)}{x-0}
=\frac{\frac{1}{n+1}+\frac{\tau}{2}}{\frac{1}{n+1}+\tau}
=\frac{1+\frac{\tau}{2}(n+1)}{1+\tau(n+1)}
\gl\bgar{ll} \le 1, \\[6pt] \ge \frac{1}{1+\frac{1}{n}}. \ear\r.$$
From this and the fact that $f$ is an odd function on $(-1,1)$ we
see  that the derivative $f'(0)=1\ne0$ exists. But $f^{-1}$ is
not continuous at $y_0=f(0)=0$ since $f^{-1}(A)=[-2,-1]\;\cup\;[1,2]$ and
since $0$ is an accumulation point of $A$. In particular, $f^{-1}$
cannot be differentiable at $y_0$. 

However, the restriction $f|_J:J\abb I\mi A$ of $f$ to the interval
$J=(-1,1)$ does have an inverse which is continuous at $y_0=0$ (and 0 is
indeed an accumulation point of $f(J)$, just as the proof of Theorem
\ref{Thm:DerivInverse} predicts). So by Theorem \ref{Thm:DerivInverse} 
the inverse $g:=(f|_J)^{-1}$ is differentiable at $0$ with $g'(0)=1$. 
\eex

On the other hand, if one assumes strict monotonicity (which in general is
stronger than injectivity), the continuity condition in (ii)' can be
skipped \cite[Satz 4.1.4 (vi)]{Behrends}. In other words, Theorem
\ref{Thm:DerivInverse} remains valid if we replace (ii) by the
condition 
\begin{itemize}
\item[(ii)''] $f$ is strictly monotonic on $I$. 
\end{itemize}

This condition is satisfied for the restriction $f|_{(-1,1)}$ of the
function $f$ from Example \ref{CounterexInverse}, but not for
$f:[-2,2]\abb\re$ itself. 

{\bf Acknowledgment.} We are very grateful to Professor Stephan
Ruscheweyh and to Professor Lawrence Zalcman for several helpful
comments.

\bibliographystyle{amsplain}

\vspace{12pt}

\parbox{90mm}{
J\"urgen Grahl \\
University of W\"urzburg \\
97074 W\"urzburg\\
Germany\\
e-mail: grahl@mathematik.uni-wuerzburg.de}
\hfill\parbox{90mm}{
Shahar Nevo \\
Bar-Ilan University\\
Ramat-Gan 52900\\
Israel\\
e-mail: nevosh@math.biu.ac.il
}

\end{document}